\documentclass[fleqn,a4paper,12pt]{article}
\usepackage{amsmath,amsfonts,calrsfs,amssymb,color}

\newtheorem{Theorem}{Theorem}[section]
\newtheorem{Definition}[Theorem]{Definition}
\newtheorem{Proposition}[Theorem]{Proposition}
\newtheorem{Lemma}[Theorem]{Lemma}
\newtheorem{Corollary}[Theorem]{Corollary}
\newtheorem{Remark}[Theorem]{Remark}

\makeatletter
\@addtoreset{equation}{section}

\makeatother

\def\R{\mathbb R}
\def\N{\mathbb N}

\def\E{\mathbb E}
\def\P{\mathbb P}
\def\ds{\displaystyle}

\newcommand{\one}{1\!\!\!\;\mathrm{l}}

\title{\bf Asymptotic behavior of stochastic PDEs with random coefficients}
\author{Giuseppe Da Prato
\\Scuola Normale Superiore di Pisa,\\  Piazza dei Cavalieri
7, 56126 Pisa, Italy
\and
Arnaud Debussche\\
Ecole Normale Sup\'erieure de Cachan,\\
 antenne de
Bretagne,\\
 Campus de Ker Lann, 35170 Bruz, France
}
\date{ }

\begin{document}
\maketitle

\begin{abstract} 
We study the  long time behavior of the solution
of a stochastic PDEs with random coefficients
assuming that randomness arises in a different independent scale.   We apply the obtained results to $2D$- Navier--Stokes equations. 
\end{abstract}
\bigskip

\noindent {\bf 2000 Mathematics Subject Classification AMS}:  76D05, 60H15, 37A25, 37L55

\bigskip 

\section{Introduction and setting of the problem}

In this work, we consider partial differential equations driven by two different sources of 
randomness. One is of white noise type and the other one is smoother. This kind of problem is 
very natural. For instance, if a physical system is submitted to external random forces and if
these have different time scales the fast ones can often be approximated by a white noise. 

We thus study the following stochastic differential equation,
\begin{equation}
\label{e1.1}
\left\{\begin{array}{l}
dX=(AX+b(X)+g(X,\overline Y)dt+\sigma(X,\overline Y))dW(t),\\
\\
X(s)=x\in H,\quad s\le t,
\end{array}\right. 
\end{equation}
where  $A:D(A)\subset H\to H$ is the infinitesimal generator of a strongly continuous semigroup
$e^{tA}$,  $b:D(b)\subset H\to H$, $g:D(g)\subset H\times K\to H$,  $\sigma:  H\times K\to L_2(H)$ are suitable nonlinear mappings.

The two sources of randomness are the Wiener process $W$ and the process $\overline Y$. We assume that they are independent. More precisely, 
we are given two separable Hilbert spaces $H$ and $K$ and two filtered probability spaces $(\Omega, \mathcal F,\P,(\mathcal F_t)_{t\ge 0})$  and 
$(\widetilde{\Omega}, \widetilde{\mathcal F}, \widetilde{\P}, (\widetilde{\mathcal F}_t)_{t\ge 0})$. We shall set  $\E=\E^\P,\;{\widetilde\E}=\widetilde\E^{\widetilde\P}$.
Then $W(t)$ is a cylindrical  Wiener process on $H$ associated to the stochastic basis  $(\Omega, \mathcal F,\P,(\mathcal F_t)_{t\ge 0})$ and $\overline Y(t)$ is a $K$-valued  Markov stationary  process associated to the stochastic basis   $(\widetilde{\Omega}, \widetilde{\mathcal F}, \widetilde{\P}, (\widetilde{\mathcal F}_t)_{t\ge 0})$ independent of $W$.

Several problems   can be written as equation \eqref{e1.1}. For instance, it may describe 
the evolution of a fluid and equation \eqref{e1.1} is then an abstract form of the Navier-Stokes equation. 
Other models are reaction-diffusion equations, Ginzburg-Landau equations and so  on.
 
 In \cite{DD}, \cite{DaPrato-Roeckner}, the case when $\overline Y$ is deterministic and periodic in time
 has been studied. Since  $X$ is not an homogeneous Markov process, the notion of invariant measure does not make sense anymore. Instead, the longtime behaviour is described by an evolutionary system of measures 
$(\mu_t)_{t\in\R}$. It is periodic and is such that if the law of $X$ at time $s$ is $\mu_s$ then it is $\mu_t$
at time $t$. In fact, this evolutionary system can be constructed by disintegration of an invariant 
measure of an enlarged system which is Markov. Moreover, uniqueness, ergodicity  and mixing 
properties have been generalized to this context.

 Our aim in this 
 work is to generalize the results obtained in \cite{DD} to more general driving forces. We also 
 define systems of measures which generalize the concept of invariant measures and describe 
 the longtime behaviour. It is also obtained by disintegration of an invariant measure but in a more
 complicated way. 
 
We assume that  \eqref{e1.1}  has a unique continuous   solution. This is the case in the examples described above if the function $g(X,\overline Y)$ is Lipschitz in $X$, has polynomial growth in 
$\overline Y$ and $\overline Y$ has moments for instance.  We denote the solution by
$X(t,s,x)$. 

We also assume  that there is a continuous Markov process $Y(t,s,y),\;t\ge s,\;y\in K$ on 
 $(\widetilde{\Omega}, \widetilde{\mathcal F}, \widetilde{\P}, (\widetilde{\mathcal F}_t)_{t\ge 0})$ such that $Y(s,s,y)=y$ and
$$
Y(t,s,\overline Y(s))=\overline Y(t),\quad t\ge s.
$$

A typical example is when $Y$ is the solution of a stochastic partial differential 
equation driven by another Wiener process and $\overline Y$ is a stationary solution 
of this equation. The simplest case is given by an Ornstein Uhlenbeck process:
\begin{equation}
\label{e1.2}
Y(t,s,y)=e^{(t-s)B}y+\int_s^te^{(t-r)B}dV(r){,\quad t>s,}
\end{equation}
where $B:D(B)\subset K\to K$   is self-adjoint, strictly negative and such that $B^{-1}$ is of trace class and $V$ is a cylindrical Wiener process in $\R$ with values in $K$ independent of $W$.

 We set
$$
P_{s,t}\varphi(x)=P^{\widetilde\omega}_{s,t}\varphi(x)=\E[\varphi(X(t,s,x))],\quad\varphi\in B_b(H). 
$$
Obviously,  $P_{s,t}\varphi(x)$ is a random variable
in $(\widetilde{\Omega}, \widetilde{\mathcal F}, \widetilde{\P})$ and   for each $\widetilde\omega\in \widetilde{\Omega}$ it holds
$$
P_{s,r}P_{r,t}=P_{s,t},\quad s\le r\le t.
$$
One should not confuse $P_{s,t}\varphi(x)$ with
$\widetilde\E[P_{s,t}\varphi(x)]$. The latter 
 does not fulfil the cocycle law  since
 $X(t,s,x)$ is not a Markov process in general.

In this paper we are interested in the asymptotic behavior of $P_{s,t}\varphi(x)$. We shall proceed as follows. First we construct an enlarged homogeneous Markov process $Z(t,s,x,h)$ with state space
$H\times\mathcal K$, where $\mathcal K:=C((-\infty,0];K)$.

Then assuming that $Z$ possesses an invariant measure $\nu(dx,dh)$ we show that a disintegration of $\nu$ produces a family $(\mu_t)_{t\in\R}$ of random  probability measures on $H$ such that $\widetilde\omega$ a.s.
\begin{equation}
\label{e1.3d}
\int_HP^{\widetilde\omega}_{s,t}\varphi(x)\mu_s(dx)=\int_H\varphi(x)\mu_t(dx),\quad \varphi\in B_b(H),\; t>s.
\end{equation}
Such a family of measures is called an evolutionary system of measures. 

Then, we give sufficient conditions for uniqueness of an evolutionary system of measures. We 
generalize the classical criterion based on irreducibilty and strong Feller property. We also show
that the recent method developed in \cite{HM}   generalizes to our context.
Finally, we illustrate our results on the two dimensional Navier-Stokes.

\section{Evolutionary systems of meausres}
\label{s2}
\subsection{Construction of  a Markov process}
\label{s2.1}

  Let us first define a new Markov process on $\mathcal K:=C((-\infty,0];K)$, setting
  \begin{equation}
\label{e1.3}
H(t,s,h)(\theta)=\left\{\begin{array}{l}
Y(t+\theta,s,h(0)),\quad\mbox{\rm if}\;t+\theta\ge s,\\
h(\theta-s+t),\quad\mbox{\rm if}\;t+\theta< s.
\end{array}\right., \quad h\in \mathcal K.\\
\end{equation}
It is not difficult to check that 
\begin{equation}
\label{e1.4b}
H(t,s,h)=H(t,r,H(r,s,h)),\quad t\ge r\ge s.
\end{equation}
It follows easily that $H$ is a Markov process which is clearly homogeneous. 

We assume that for any $h\in \mathcal K$ the following equation has a unique continuous solution
\begin{equation}
\label{e1.4}
\left\{\begin{array}{lll}
dX^h&=&(AX^h+b(X^h)+g(X^h,H(t,s,h)(0))dt\\
\\
&&+\sigma(X^h,H(t,s,h)(0)))dW(t),\\
\\
X^h(s)&=&x\in H,\quad s\le t.
\end{array}\right. 
\end{equation}

Then we define the spaces
$$
( \Omega_1, \mathcal F_1, \P_1, (\mathcal F_{1,t})_{t\ge 0})=(\Omega\times\widetilde{\Omega}, \mathcal F\times\widetilde{\mathcal F}, \P\times\widetilde{\P}, (\mathcal F_t\times\widetilde{\mathcal F}_t)_{t\ge 0}),
$$
(we shall denote by $\E_1$ the expectation in this space) and
$$
\mathcal H=H\times \mathcal K,
$$
and we consider  the homogeneous Markov process
$$
Z(t,s,x,h)=(X^h(t,s,x),H(t,s,h)).
$$
We denote by $Q_t,\;t\ge 0,$ and $R_t, \; t\ge 0$ the  transition semigroups associated to $Z$ and $H$ respectively.

We also denote by $P^h_{t,s}$ the transition operators associated to $X^h$.

\subsection{Stationary processes}
\label{s2.2}

If 
$$
h(\theta)=\overline Y(s+\theta),\quad \theta\le 0,
$$
then
$$
H(t,s,h)=Y(t,s,h(0))=\overline Y(t+\cdot),\quad t\ge 0.
$$
Therefore $H(t,s,h)(0)=\overline Y(t),$
so that
$$
X^h(t,s,x)=X(t,s,x),
$$
where $X(t,s,x)$ is the solution to \eqref{e1.1}. Moreover
$H(t,s,h)(\theta)=\overline Y(t+\theta)$ so that 
\begin{equation}
\label{e1.5b}
H(t,s,h) =\tau_t \overline Y|_{(-\infty,t]},
\end{equation}
where, for a function $f$ defined on $(-\infty,t]$, $\tau_tf$ is the function defined on $(-\infty,0]$
by $\tau_t f(\theta)=f(t+\theta)$. It follows that
$$
Z(t,s,x,h)=(X(t,s,x),\tau_t \overline Y|_{(-\infty,t]}).
$$

Note also that the stochastic process $H(t,s,h)$ is stationary since the laws of 
$\tau_t \overline Y|_{(-\infty,t]}$ and $\tau_s \overline Y|_{(-\infty,s]}$ coincide for any $t,s$. 
Moreover, it does not depend on $s$. We shall denote it by $\overline H$, and by $\lambda$  its invariant law.

Conversely, let $\overline H$ be a stationary process associated to $H$. 
Then 
$$
\overline H(t)(\theta ) = Y(t+\theta,s ,\overline H (s)(0))
$$
for all $t+\theta \ge s$. Therefore, if we define
$$
\overline Y(t)=\overline H(t)(0)
$$
we have $\overline Y(t) = H(\tilde t)(\tilde \theta)$ provided $\tilde t+\tilde \theta = t$. In particular,
$\overline Y(t)=\overline H(0)(t)$ for $t\le 0$.

Note also that  $\overline Y(t) = \overline H(t)(0) = Y(t,s,\overline H(s)(0)) =Y(t,s, \overline Y(s))$. Since $\overline Y$ is clearly stationary, it follows that it is a stationary process associated to $Y$.

This shows that we have  a one to one correspondance between stationary processes for $H$ and $Y$.

Let $\lambda$ be an invariant measure for $H$:
$$
\int_{\mathcal K}\alpha(H(t,s,h))\lambda(dh)=\int_{\mathcal K}\alpha(h)\lambda(dh),\; \alpha\in C_b(\mathcal K).
$$
 Take $\alpha(h)=\gamma(h(0))$, $h\in \mathcal K,
$
where $\gamma\in C_b(K)$.
Since $H(t,0,h)(0)=Y(t,0,h(0))$,  we deduce
\begin{equation}
\label{e1.10b}
\int_  K \widetilde{\E}[\gamma(Y(t,0,h(0)))]\lambda(dh)=\int_  K\gamma(h(0))\lambda(dh),\quad\forall\;\gamma\in C_b(K).
\end{equation}

If we denote by $\lambda_0$ the image measure of $\lambda$ by the mapping
$$
\mathcal K\to K,\;h\mapsto h(0),
$$
we deduce from \eqref{e1.10b} that
\begin{equation}
\label{e1.11b}
\int_\mathcal K \widetilde{\E}[\gamma(Y(t,0,\xi))]\lambda_0(d\xi)=\int_\mathcal K\gamma(\xi)\lambda_0(d\xi),\quad\forall\;\gamma\in C_b(K).
\end{equation}
Therefore $\lambda_0$ is an invariant measure for $Y$.

We now prove that ergodicity is transferred to $\overline H$. In the proof of this result, we also prove that 
the law of the process $\overline Y$ determines the law of $\overline H$.

\begin{Proposition}
\label{p5.1i}
Assume that $Y$ has a unique invariant measure, then the process   $\overline H$ is stationary and ergodic.
\end{Proposition}
{\bf Proof}. Let $\lambda$ be an invariant law of $\overline H $. Assume that
$$
\lambda=\eta\lambda^1+(1-\eta)\lambda^2,
$$
where $\eta\in[0,1]$ and $\lambda^1,\lambda^2$ are invariant for $\overline H$. Given $\theta\in \R$ consider the mapping
$$
\mathcal K\to K,\quad h\mapsto h(\theta)
$$
and denote by $\lambda_\theta,\;\lambda^1_\theta,\;\lambda^2_\theta$ the image measures of $\lambda,\lambda^1,\lambda^2$ respectively  by this mapping.
By the discussion above, $\lambda_0,\;  \lambda_0^1,\; \lambda_0^2$ are invariant for $Y$. Therefore
$$
\lambda_0=\lambda_0^1=\lambda_0^2.
$$
We have in addition for $\gamma_1,\gamma_2\in C_b( K)$ and $\theta_1\in \R,\;\theta_2\ge 0$
$$
\begin{array}{l}
\ds \int_{\mathcal K}\gamma_1(h(\theta_1))\gamma_2
(h(\theta_1+\theta_2))\lambda(dh)
\\
\\
\ds =\int_{\mathcal K}\widetilde\E[\gamma_1(Y(t+\theta_1,0,h(0)))\gamma_2
(Y(t+\theta_1+\theta_2,0,h(0)))]\lambda(dh)\\
\\
\ds=\int_{K}\widetilde\E[\gamma_1(Y(t+\theta_1,0,y))\gamma_2
(Y(t+\theta_1+\theta_2,0,y))]\lambda_0(dy)\\
\\
\ds=\int_{K}\widetilde\E[\gamma_1(y)\gamma_2
(Y(\theta_2,0,y))]\lambda_0(dy).
\end{array} 
$$
We have used the fact that $\lambda_0$ is the image law of $\lambda$  and the invariance of $\lambda_0$.

Therefore the joint law of $h(\theta_1),h(\theta_1+\theta_2)$ is given. We argue similarly for $\lambda^1$ and $\lambda^2$ and for $\theta_1,\theta_1+\theta_2,..., \theta_1+\cdots +\theta_n$ and deduce that $\lambda,\lambda^1,\lambda^2$ have the same images laws for applications
$$
\mathcal K\to K^n,\quad h\mapsto (h(\theta_1),h(\theta_1+\theta_2),..., h(\theta_1+\cdots +\theta_n)).
$$
This yields $\lambda=\lambda^1=\lambda^2$ and so, $\lambda$ is an extremal point and an ergodic invariant measure. $\Box$

\subsection{Evolutionary systems of measures}
\label{s2.3}

We now define the main new object of our article. We expect that the long time behaviour of the solutions
\eqref{e1.1} is described by a random family of measures $(\mu_t)_{t\in\R}$ satisfying 
$\widetilde\P$ a.s.
\begin{equation}
\label{e3.1f}
\int_H\E[\varphi(X(t,s,x))]\mu_{s}(dx)=\int_H\varphi(x)\mu_{t}(dx),
\end{equation}
for all $\varphi\in C_b(H)$, $t>s$. This is a natural generalisation of the notion of invariant
measure. 

This definition is too general and we need to make some restrictions on the system $(\mu_t)_{t\in\R,\widetilde{\omega}\in \widetilde{\Omega}}$. Without loss of generality, we assume that the forcing is associated to 
a shift $(\tau_t)_{t\in\R}$: 
$\overline Y(t,\widetilde\omega) = \overline Y(s,\tau_{t-s}\widetilde\omega)$,
$t\ge s,\; \widetilde\omega\in\widetilde\Omega$. 
We consider the stationary process $(\overline H(t))_{t\in\R}$
in $\mathcal K$  constructed above and defined by $\overline H (t)=\tau_t \overline Y|_{(-\infty,t]}$
and denote by $\lambda $ its invariant measure.
Then clearly $\overline H(t,\widetilde\omega) = \overline H(s,\tau_{t-s}\widetilde\omega)$,
$t\ge s,\; \widetilde\omega\in\widetilde\Omega$.

We assume that $(\mu_t)_{t\in\R,\widetilde{\omega}\in \widetilde{\Omega}}$ satisfies 
the following consistency assumption:
\begin{equation}
\label{e3.1g}
\mu_t(\widetilde\omega)=\mu_s(\tau_{t-s}\widetilde\omega),\; \mbox{a.s. }t\ge s,\; 
\mbox{and } \widetilde \P\mbox{ a.s. }\widetilde\omega\in\widetilde\Omega.
\end{equation}
This says that two solutions of \eqref{e1.1} with the same past define the same evolution. 

It is also natural to assume that, for each $t\in\R$,  $\mu_t$ is measurable with respect to the 
$\sigma$-field generated by $\overline Y(s),\; s\le t$. 

A random family of measures $(\mu_t)_{t\in\R}$ satisfying these properties is called an {\em evolutionary system of measures}.

\section{Invariant measures for $Z$ are associated to evolutionary systems of measures}
\label{s3}

Let us assume that $Z$ possesses an invariant measure $\nu$:
\begin{equation}
\label{e1.5}
\begin{array}{l}
\ds\int_{\mathcal H}\E_1[\Phi(X^h(t,s,x),H(t,s,h))]\nu(dx,dh)=\int_{\mathcal H} \Phi(x,h)\nu(dx,dh) ,
\; \Phi\in C_b(\mathcal H).
\end{array}
\end{equation}
Let us consider a  disintegration of $\nu$
$$
\nu(dx,dh)=\mu_h(dx)\lambda(dh).
$$ 
\begin{Lemma}
\label{l1.2}
$\lambda$ is an invariant measure for $H$, that is
\begin{equation}
\label{e1.8b}
\int_{\mathcal K}\widetilde{\E}[\alpha(H(t,0,h))]\lambda(dh)=\int_\mathcal K\alpha(h)\lambda(dh),\quad\forall\;  t\ge 0,\;\alpha\in C_b(\mathcal K).
\end{equation}

\end{Lemma}
{\bf Proof}.
Let   $\alpha\in C_b(\mathcal K)$. Then by \eqref{e1.5} with $s=0$ we have
$$
\begin{array}{lll}
\ds\int_{\mathcal H}\E_1[\alpha(H(t,0,h))]\nu(dx,dh)&=&\ds\int_{\mathcal H}\alpha(h)\nu(dx,dh)\\
\\
&=&\ds\int_\mathcal K\alpha(h)\lambda(dh).
\end{array} 
$$
But
$$
\int_{\mathcal H}\E_1[\alpha(H(t,0,h))]\nu(dx,dh)=\int_{\mathcal K}\widetilde{\E}[\alpha(H(t,0,h))]\lambda(dh).
$$
Therefore $\lambda$ is   invariant   for $H$. $\Box$

Let  $\overline H$ be a stationary process associated to $H$ with invariant law  $\lambda$.
We set 
$$
\mu_t=\mu_{\overline H(t)}, \; t\in \R.
$$
\begin{Theorem}
\label{t3.1f}
The family $(\mu_t)_{t\in\R}$ is an evolutionary 
system of measures.
\end{Theorem}
{\bf Proof}.  We prove that \eqref{e3.1f} holds. The other properties of evolutionary systems are clearly satisfied.

Choosing for $\varphi \in C_b(H)$, $\alpha\in C_b(\mathcal K)$:
$$
\Phi(x,h)=\varphi(x)\alpha(h),
$$
we have by \eqref{e1.5}
$$
\begin{array}{l}
\ds\int_{\mathcal H}\E_1[\varphi(X^h(t,s,x))\alpha(H(t,s,h)))]\mu_h(dx)\lambda(dh)\ds=\int_{\mathcal H} \varphi(x)\alpha(h)\mu_h(dx)\lambda(dh).
\end{array}
$$
Moreover, by the invariance of $\lambda$ for the process $H$ we get 
$$
\begin{array}{l}
\ds\int_{\mathcal H} \varphi(x)\alpha(h)\mu_h(dx)\lambda(dh)\ds=\int_{\mathcal H} \widetilde{\E}[\varphi(x)\alpha(H(t,s,h))\mu_{H(t,s,h)}](dx)\lambda(dh).
\end{array}
$$
This yields 
\begin{equation}
\label{e2.2b}
\begin{array}{l}
\ds\int_{\mathcal H}\widetilde{\E}\left\{\E[\varphi(X^h(t,s,x))\alpha(H(t,s,h)))]\right\}\mu_h(dx)\lambda(dh)\\
\\
\ds=
\int_{\mathcal H} \widetilde{\E}[\varphi(x)\alpha(H(t,s,h))\mu_{H(t,s,h)}(dx)\lambda(dh).
\end{array} 
\end{equation}
Since $\overline H(s)$ has law $\lambda$, $H(t,s,\overline H(s))=\overline H(t)$ and $X^{\overline H(s)} (t,s,x)=X(t,s,x)$, we can rewrite \eqref{e2.2b} as
$$
\begin{array}{l}
\ds\int_{H} \widetilde{\E}\left\{\E[\varphi(X(t,s,x))\alpha(\overline H(t))]\right\}\mu_{\overline H(s)}(dx)=
\int_{H} \widetilde{\E}[\varphi(x)\alpha(\overline H(t))\mu_{\overline H(t)}(dx).
\end{array}
$$
Note that $\int_H\E[\varphi(X(t,s,x))]\mu_{\overline H(s)}(dx)$ and $\int_H\varphi(x)\mu_{\overline H(t)}(dx)$ are $\overline H(t)$-measurable (their are functions of $\overline H(t)$) and we deduce by the arbitrariness of $\alpha$  that
$$
\int_H\E[\varphi(X(t,s,x))]\mu_{\overline H(s)}(dx)=\int_H\varphi(x)\mu_{\overline H(t)}(dx),\quad\widetilde\P\mbox{\rm -a.s.},
$$
which yields \eqref{e3.1f} since $\mu_t:=\mu_{\overline H(t)}$, $t\in\R$.

However, the set of all $\widetilde{\omega}\in \widetilde{\Omega}$ for which identity \eqref{e3.1f} holds depends of $t,s,\varphi$. Modifying the disintegration of $\nu$, we can easily get rid of this dependence. We proceed as in \cite{DD}. First, taking $\varphi=\one_C$ for $C$ in a countable 
set generating Borel sets of $H$, we can show that \eqref{e3.1f}
holds in fact for every $\varphi\in B_b(H)$ for almost every $\widetilde\omega$ depending only on $t,s$. From  Fubini's theorem we find that
$\widetilde \P$-a.s.
\begin{equation}
\label{e2.9b}
\int_H\E[\varphi(X(t,s,x))]\mu_{s}(dx)=\int_H\varphi(x)\mu_{t}(dx),\quad\mbox{\rm for almost all $s\le t$}.
\end{equation}
Choose $s_n\downarrow-\infty$ such that
$$
P^*_{s_n,t}\mu_{s_n}=\mu_t,\quad\mbox{\rm for almost all $t\ge s_n$} 
$$
and set
$$
\widetilde{\mu}_t^n=P^*_{s_n,t}\mu_{s_n},\quad t\ge s_n.
$$ 
From the continuity of $t\mapsto X(t,s,x)$ for almost every $\widetilde{\omega}$, we deduce that 
$t\mapsto \widetilde{\mu}_t^n$ is continuous. Moreover
$
\widetilde{\mu}_t^n=\widetilde{\mu}_t^{n+1},\quad\mbox{\rm a.s.}\;t\ge s_n,
$
so that by continuity
$$
\widetilde{\mu}_t^n=\widetilde{\mu}_t^{n+1},\quad\forall\;t\ge s_n.
$$
Now we define
$$
\widetilde{\mu}_t=\widetilde{\mu}_t^n,\quad\forall\;t\ge s_n.
$$
Obviously, $\widetilde{\mu}_t=\mu_t$ for almost all $t\in\R$,
so that
$$
P^*_{s_n,t}\widetilde{\mu}_{s_n}=\widetilde{\mu}_t,\quad\mbox{\rm a.s.}\;t\ge s_n.
$$
By the continuity in $t$ we deduce for all $t\in\R$
$$
P^*_{s_n,t}\widetilde{\mu}_{s_n}=\widetilde{\mu}_t,\quad\forall\;t\ge s_n
$$
and for $t\ge s\ge s_n$ we have
$$
P^*_{s,t}\widetilde{\mu}_{s}=P^*_{s,t}(P^*_{s_n,s}\widetilde{\mu}_{s_n})=P^*_{s_n,t}\widetilde{\mu}_{s_n}=\widetilde{\mu}_t.
$$
Therefore we can conclude that $\widetilde \P$-a.s.
$$
P^*_{s,t}\widetilde{\mu}_{s}=\widetilde{\mu}_t,\quad\forall\;  t\ge s,
$$
as claimed. Since $\widetilde{\mu}_t=\mu_t$ for almost all $t\in\R$, the consistency relation is still 
satisfied by $\widetilde \mu_t$.
$\Box$

Let us now see how an evolutionary system of measures $(\widetilde\mu_t)_{t\in\R,\widetilde{\omega}\in \widetilde{\Omega}}$ yields an invariant measure for $Z$.
By definition, $\widetilde\mu_t$ is measurable  with respect to the 
$\sigma$-field generated by $\overline Y(s),\; s\le t$. This implies 
that there exists  a measurable function
$h\mapsto \mu^t_h$ such that $\widetilde\mu_t=\mu^t_{\overline H(t)}$. By the consistency
assumption \eqref{e3.1g}, we have $\widetilde\P$-a.s.  for $t\ge s$ in a set $I$ of full measure,
$$
\mu^t_{\overline H(t,\widetilde\omega)}= \widetilde\mu_t(\widetilde\omega)
=\widetilde\mu_s(\tau_{t-s}\widetilde\omega)=\mu^s_{\overline H(s,\tau_{t-s}\widetilde\omega)}
=\mu^s_{\overline H(t,\widetilde\omega)}
$$
so that $\mu^t$ does not depend on $t$ for $t\in I$. We choose $s_0\in I$, set $\mu_h=\mu^{s_0}_h$ and define 
$$
\nu(dx,dh)=\mu_h(dx)\lambda(dh).
$$
Then, since $\mathcal L(\overline Y(\cdot,s_0))=\lambda$ and $X^{\overline H(s_0)}(t,s_0,x)
=X(t,s_0,x)$, we have for $\psi\in B_b(\mathcal H)$, $t\in I$,
$$
\begin{array}{ll}
\ds\int_{\mathcal H} Q_{t-s_0} \psi(x,h)\nu(dx,dh)&\ds=\int_{\mathcal H}\E_1[\psi(X^h(t,s_0,x)),H^h(t,s_0,h))]\mu_h(dx)\lambda(dh)\\
\\
&\ds=\int_{H}\E_1[\psi(X(t,s_0,x)),\overline H(t))]\mu_{\overline H(s_0)}(dx).
\end{array} 
$$
Therefore, for $\psi$ of the form $\psi(x,h)=\varphi(x)\alpha(h)$,
$$
\begin{array}{ll}
\ds\int_{\mathcal H} Q_{t-s_0} \psi(x,h)\nu(dx,dh)
&\ds=\int_{H}\E_1[\varphi(X(t,s_0,x))\alpha(\overline H(t))]\widetilde\mu_{s_0}(dx)\\
\\
&\ds= \widetilde\E\left[\int_H\E(\varphi(X(t,s_0,x)))\widetilde\mu_{s_0}(dx) \alpha(\overline H(t))  \right]
\\
\\
&\ds= \widetilde\E \int_H\varphi(x)\widetilde\mu_{t}(dx) \alpha(\overline H(t)) \\
\\
&\ds= \widetilde\E \int_H\varphi(x)\mu_{\overline H(t)}(dx) \alpha(\overline H(t))\\
\\
&\ds=\int_{\mathcal H}\varphi(x)\alpha(h) \mu_h(dx)\lambda(dh).\\
\\
 &\ds=\int_{\mathcal H}\psi(x,h)\nu(dx,dh). 
\end{array} 
$$
Since the left hand side is a continuous function of $t$, we deduce that the equality holds for all $t\in \R$ and that $\nu$ is an invariant measure for $Z$.

It follows from our discussion that the correspondance $\nu \mapsto ((\mu_t)_{t\in\R},\lambda)$ is a bijection.
In particular, if there exists a unique invariant measure with marginal $\lambda$ for $Z$, there exists a unique evolutionary system of measure.

\section{A simple example}
\label{s2.4}

We consider the following special form of \eqref{e1.1} in $\R$.
\begin{equation}
\label{e5.1h}
\left\{\begin{array}{lll}
dX&=&(-X +\overline Y)dt+ dW, \; t\ge s \\
\\
X(s)&=&x\in H,
\end{array}\right. 
\end{equation}
where  $\overline Y$
is the stationary process
\begin{equation}
\label{e5.2h}
\overline Y(t)=\int_{-\infty}^te^{-(t-r)}dV(r), {\quad t\in\R}
\end{equation}
 and $V$, $W$ are independent Wiener processes.
 
 We consider the homogeneous Markov process introduced above
$$Z(t,s,x,h)=(X^h(t,s,x),H(t,s,h)),\quad t\ge s,\;x\in H,\; h\in \mathcal K.$$

Let us write explicitly  $Z(t,s,x,h).$
We have
\begin{equation}
\label{e5.3h}
X^h(t,s,x)=e^{-(t-s)}x+\int_s^te^{-(t-r)}Y(r,s,h(0))dr + \int_s^t e^{-(t-r)}dW(r).
\end{equation}
Then we have
\begin{equation}
\label{e4}
\left\{\begin{array}{l}
\ds\lim_{s\to-\infty}X^h(t,s,x)=\int_{-\infty}^te^{(t-r)A}\overline Y(r)dr+ \int_{-\infty}^t e^{(t-r)A}dW(r):=\zeta_t\\
\\
\ds \lim_{s\to-\infty}H(t,s,h)(\theta)=\bar Y(t+\theta).
\end{array}\right.
\end{equation}
This implies that  the probability measure on $H\times\mathcal K$,
$$
\nu=\mathcal L(\zeta_0,\overline H(0)),
$$
is the unique invariant measure for $Z$. In other words we have
$$
\int_{H\times\mathcal K}\Phi(x,h)\nu(dx,dh)=\E_1[\Phi(\zeta_0,\overline H(0))].
$$
Let $\nu(dx,dh)=\mu_x(dx)\lambda(dh)$ be a disintegration of $\nu$. Then we have
$$
\lambda=\mathcal L(\overline H(0))
$$
and
$$\mu_t= \mathcal N\left(\int_{-\infty}^0e^{-\theta}\overline Y(t+\theta)d\theta, \frac12\right).$$
We see on this simple example that it is necessary to parametrize the evolutionary system 
of measure by the whole history of the driving process.

\section{Ergodicity in the regular case}
\label{s4}

We assume for simplicity in this section that the Markov process $Y$ has a unique invariant measure. It is then the invariant law of $\overline Y$ and $\overline Y$ is ergodic. We construct the ergodic
process $\overline H$ as above and denote by $\lambda$ its invariant law.

We generalize the famous Doob criterion (see for instance \cite{DPZ2}) of ergodicity to evolutionary system of measures.

Let us set
$$
\pi_{s,t}(x,E)=\pi^{\widetilde\omega}_{s,t}(x,E)=P_{s,t}\one_E(x),\quad\forall\;I\in \mathcal B(H).
$$
We say that $P_{s,t}$ is {\em regular} at $\widetilde\omega$ if  for any $s<t$ we have
$$
\pi^{\widetilde\omega}_{s,t}(x,\cdot)\sim \pi^{\widetilde\omega}_{s,t}(y,\cdot)\quad\mbox{\rm for all}\;x,y\in H.
$$
We notice the following straightforward identity,
\begin{equation}
\label{e4.1g}
\pi_{s,t}(x,E)=\int_H\pi_{s,s+h}(x,dy)\pi_{s+h,t}(y,E),\quad s+h<t,\;h>0,\;E\in \mathcal B(H). 
\end{equation}
We say that  $P_{s,t}$ is strong Feller at $\widetilde\omega$ if for each $s<t$, $P_{s,t}$ maps $\mathcal B(H)$
to $C_b(H)$. It is irreducible at $\widetilde\omega$ if,  for any $s<t$, $x\in H$ and $O\subset H$ open, $\pi_{s,t}(x_0,O)>0$. The proofs of the following Propositions \ref{p4.1g} and \ref{p4.2g} are completely similar to that of \cite[Proposition 4.1, Proposition 4.2]{DD}.
\begin{Proposition}
\label{p4.1g}
If the transition semigroup $P_{s,t}$ is strong Feller and irreducible at $\widetilde\omega$, then it is regular at $\widetilde\omega$.
\end{Proposition}
\begin{Proposition}
\label{p4.2g}
 Assume that $P_{s,t}$ is regular at $\widetilde\omega$ and that it possesses an invariant set  of probabilities $\mu_t,\;t\in\R$. Then $\mu_t (\widetilde\omega)$ is equivalent to $\pi^{\widetilde\omega}_{s,t}(x,\cdot)$ for all $s<t$ and 
 $x\in H$.
\end{Proposition}

Next we prove the following theorem which can be used in several applications provided the
noise $W$ is non degenerate.

\begin{Theorem}
\label{t4.3g}
We assume that $P_{s,t}$ is regular for almost all $\widetilde\omega$ and that the Markov process $Y$ has a unique invariant measure.
Then $Q_t$ has at most one invariant measure which is in addition ergodic. It follows that there exists a unique evolutionary system
of measures.
\end{Theorem}
{\bf Proof}. It suffices to prove that all invariant measures of $Q_t$ are ergodic. Let $\nu(dx,dh)=
\mu_h(dx)\lambda(dh) $ be an invariant measure for $Q_t$
and $\Gamma\in \mathcal B(\mathcal H)$  be an invariant set for $\nu$:
$$
Q_t\one_\Gamma(x,h)=\one_\Gamma(x,h),\quad\nu\mbox{\rm-a.e.}
$$
Then for any $\Phi\in B_b(\mathcal H)$ we have
$$
\int_\mathcal H\left(Q_t\one_\Gamma\right)(x,h)\Phi(x,h)\nu(dx,dh)=\int_\mathcal H\one_\Gamma(x,h)\Phi(x,h)\nu(dx,dh).
$$
Setting
$$
\Gamma_h=\{x\in H:\;(x,h)\in \Gamma\},
$$
we have $\one_\Gamma(x,h)=\one_{\Gamma_h}(x)$ and
$$
\begin{array}{l}
\ds\int_\mathcal H\E_1\left(\one_{\Gamma_{H(t,0,h)}} (X^h(t,0,x)   \right)\Phi(x,h)\nu(dx,dh)
\ds=\int_\mathcal H\one_{\Gamma_h}(x)\Phi(x,h)\nu(dx,dh).
\end{array}
$$
On the other hand,
$$
\begin{array}{l}
\ds\int_\mathcal H\E_1\left(\one_{\Gamma_{H(t,0,h)}} (X^h(t,0,x)   \right)\Phi(x,h)\nu(dx,dh)\\
\\
\ds=\int_\mathcal H\widetilde \E\left(P^h_{0,t}\one_{\Gamma_{H(t,0,h)}}    \right)(x)\Phi(x,h)\mu_h(dx)\lambda(dh)\\
\end{array}
$$
We deduce
$$
\widetilde\E\left(P_{0,t}^h\one_{\Gamma_{H(t,0,h)}}\right)=\one_{\Gamma_h}, \;  \nu\mbox{-a.s.}
$$
so that
$$
P_{0,t}^h\one_{\Gamma_{H(t,0,h)}}=\one_{\Gamma_h}, \; \widetilde \P\times \nu\mbox{-a.s.}
$$
and  
$$
P_{0,t}^{{\overline H(0)}}\one_{\Gamma_{\overline H(t)}}=\one_{\Gamma_{\overline H(0)}},\quad\mu_{\overline H(0)}\times\widetilde \P \mbox{\rm -a.s.}
$$
Therefore, $\widetilde \P$-a.s., if $\mu_{\overline H(0)}(\Gamma_{\overline H(0)})\neq 0$ then for $\mu_{\overline H(0)}$ almost every $x\in \Gamma_{\overline H(0)}$ we have $P_{0,t}^{\overline H(0)}1_{\Gamma_{\overline H(t)}}(x)=1$, equivalently
$\pi_{0,t}^{\overline H(0)}(x,\Gamma_{\overline H(t)})=1$. Since, by assumption, $P_{0,t}^{{\overline H(0)}}$ is regular, we have
$$
\pi_{0,t}^{\overline H(0)}(y,\Gamma_{\overline H(t)})=1,\quad\forall\;y\in H,
$$
and so,
$$
\one_{\Gamma_{\overline H(0)}}(y)=1,\quad\mu_{\overline H(0)}\mbox{\rm -a.s.}
$$
Therefore
$$
\mu_{\overline H(0)}(\Gamma_{\overline H(0)})=1.
$$
We have proved that, $\widetilde \P$-a.s., 
$$
\mu_{\overline H(0)}(\Gamma_{\overline H(0)})=\; \mbox{\rm $0$ or $1$}.
$$
Similarly we show that
$$
\mu_{\overline H(t)}(\Gamma_{\overline H(t)})=\; \mbox{\rm $0$ or $1$},\;\quad\forall\;t\in\R.
$$
Set
$$
\widetilde\Omega_1:=\{\widetilde\omega\in \widetilde\Omega:\;\mu_{\overline H(1)}(\Gamma_{\overline H(1)})=1\}.
$$
Then
$$
\tau_1\widetilde\Omega_1:=\{\widetilde\omega\in \widetilde\Omega:\;\mu_{\overline H(0)}(\Gamma_{\overline H(0)})=1\},
$$
where $\theta_t$ is the shift defined in sectionÊ \ref{s2.3}

Write now
$$
\begin{array}{l}
\ds\mu_{\overline H(0)}(\Gamma_{\overline H(0)})=\int_H\one_{\Gamma_{\overline H(0)}}d\mu_{\overline H(0)}=\int_HP^{\overline H(0)}_{0,1}\one_{\Gamma_{\overline H(1)}}d\mu_{\overline H(0)}
\\
\\
\ds =\int_H\one_{\Gamma_{\overline H(1)}}d\mu_{\overline H(1)}=\mu_{\overline H(1)}(\Gamma_{\overline H(1)}).
\end{array}
$$
Therefore $\tau_1\widetilde\Omega_1=\widetilde\Omega_1$ and the ergodicity implies $\widetilde\P(\widetilde\Omega_1)=0$ or $1$.

If $\widetilde\P(\widetilde\Omega_1)=0$   we have
$
\nu(\Gamma)=\widetilde\E\mu_{\overline H(0)}(\Gamma_{\overline H(0)})=0
$ and if $\widetilde\P(\widetilde\Omega_0)=1$   we have
$
\nu(\Gamma)=\widetilde\E\mu_{\overline H(t)}(\Gamma_{\overline H(t)})=1.
$
$\Box$
\begin{Remark}
\em The proof of Theorem 4.3  in \cite{DD} is not complete. In fact in that theorem we have  proved only that $\nu_h(\Gamma_h)=0$ or $1$. One has to use the same argument as before to arrive at the conclusion.
$\Box$
\end{Remark}

\section{Uniqueness by asymptotic strong Feller property}
\label{s5}

We assume again that $\overline Y$ is ergodic and that the process $Y(\cdot,s,y)$ has a unique invariant measure. 
The following definition is a natural generalization to $P_{s,t}$
of a concept introduced in \cite{HM}.
\begin{Definition}
We say that  $P_{s,t}$ is {\em asymptotic strong Feller (ASF)}
at $x\in H$ if there is a sequence of  pseudo--metrics $\{d_n\}$ 
  on $H$ such that
\begin{equation}
\label{e5.1i}
d_n(x_1,x_2)\uparrow 1,\quad \forall\; x_1\ne x_2,
\end{equation}
and a sequence $\{t_n\}$
of positive numbers such that,
\begin{equation}
\label{e5.2i}
\lim_{\gamma\to 0}\;\limsup_{n\to\infty}\widetilde\E\left(\sup_{y\in B(x,\gamma)}
W_{d_n}(\pi_{s,s+t_n}(x,\cdot),
\pi_{s,s+t_n}(y,\cdot))\right)=0.
\end{equation}
$P_{s,t}$ is called {\em asymptotically strong Feller (ASF)}, 
if it is asymptotically strong Feller at any   $x\in
H$ and $s_0\in\R$.
\end{Definition}
As in \cite{DD} section 5.1, $W_d$ denotes the Wasserstein metric on the space of probability measures on 
$H$ associated to a pseudo metric $d$. Recall that $d_n(x_1,x_2)\uparrow 1,\quad \forall\; x_1\ne x_2,$ implies that $W_{d_n}(\mu_1,\mu_2)\to \|\mu_1-\mu_2\|_{TV}$, the total variation distance between $\mu_1$ and $\mu_2$ (see \cite{HM}, Corollary 3.5).
\begin{Remark}
\em  The previous definition is independent of $s$ because by stationarity of 
$\overline Y$
$$
\widetilde\E \left[\sup_{y\in B(x,\gamma)}W_{d_n}(\pi_{s,s+t_n}(x,\cdot),\pi_{s,s+t_n}(y,\cdot))   \right],
$$
is independent of $s$.

\end{Remark}
The following result can be proved as in \cite{HM}.
\begin{Proposition}
\label{p5.4i}
Assume that for some $s>0$ there exist
  $t_n\uparrow +\infty$, $\delta_n\to 0$ and $C(s, |x|)$  locally bounded  with respect to $|x|$ and such that for $\gamma<1$,
\begin{equation}
\label{e5.5}
\widetilde\E\left(\sup_{y\in B(x,\gamma)}|DP_{s,s+t_n}\varphi(x)|\right)\le C( s, |x|)(\|\varphi\|_\infty+\delta_n\|D\varphi\|_\infty).
\end{equation}
Then $P_{s,t}$ is ASF.
\end{Proposition}
\begin{Lemma}
\label{l5.5i}
 Let $\nu^1(dx,dh)=\mu^1_h(dx)\lambda(dh)$, $\nu^2(dx,dh)=\mu^2_h(dx)\lambda(dh)$ be two invariant measures for $Q_t$ such that
$\nu_1$ and $\nu_2$ are singular.
 Then $\mu^1_h$ and $\mu^2_h$ are singular for  $\lambda$ almost all $h\in \mathcal K$.
\end{Lemma}
{\bf Proof}. Let $A,B\in \mathcal B(\mathcal K)$ such that $\nu^1(A)=\nu^2(B)=1$ and $A\cap B=\varnothing.$ For each $h\in \mathcal K$ we define
$$
A_h=\{x\in H:\;(x,h)\in A\},\;
B_h=\{x\in H:\;(x,h)\in B\}.
$$
Then  we have $A_h\cap B_h=\varnothing$ and $\mu^1_h(A_h)=\mu^2_h(B_h)=1$ because
$$
\nu^1(A)=\int_{\mathcal K}\mu^1_h(A_h)\lambda(dh)=1,\;
\nu^2(B)=\int_{\mathcal K}\mu^2_h(B_h)\lambda(dh)=1.
$$
$\Box$

Using the same proof of for  Lemma 5.6 in \cite{DD}, we prove the following result.
\begin{Lemma}
\label{l5.6i}
Let $d\le 1$ be a pseudo-metric on $\mathcal H$ and let $\nu^1$ and $\nu^2$ be two invariant measures for $Q_t$ with the same marginal 
$\lambda$. Let us denote by $(\mu^1_t)_{t\in\R}$
and $(\mu^2_t)_{t\in\R}$ the system of random measures constructed in Section \ref{s3}. Then we have
\begin{equation}
\label{e5.4i}
W_d(\mu^1_{t+s},\mu^2_{t+s})=1-\mu^1_s(A)\wedge
\mu^2_s(A)\left(1-\sup_{y,z\in A} W_d(\pi_{s,t+s}(y,\cdot)-\pi_{s,t+s}(z,\cdot))  \right).
\end{equation}
\end{Lemma}

\begin{Theorem}
\label{t5.7i}
Assume that $(P_{s,t})_{t\ge s}$ is ASF and that there is $x_0\in H$ such that all invariant measure $\nu$ of $Q_t$ 
\begin{equation}
\label{support}
x_0\in\;\mbox{\rm supp}\;\mu_h,\quad \lambda \mbox{ -a.s.}
\end{equation}
Then there exists at most one invariant measure  for $(Q_t)_{t\ge 0}$.
\end{Theorem}
{\bf Proof}. It is enough to show that two ergodic invariant measures of $(Q_t)_{t\ge 0}$, $\nu^1,\nu^2$    are identical. If not then $\nu^1$ and $\nu^2$ are necessarily singular and so, by Lemma \ref{l5.5i}, $\mu_h^1,\mu_h^2$ are singular 
$\lambda$-a.s. Since $Q_t$ is ASF there exists $\gamma_0>0$, $n\in\N$ such that
$$
\widetilde\E\left[\sup_{y,z\in B(x_0,\gamma_0)}(W_{d_n}(\pi_{0,t_n}(y,\cdot),\pi_{0,t_n}(z,\cdot))   \right]<\frac14.
$$
We deduce
$$
\widetilde\P\left[\sup_{y,z\in B(x_0,\gamma_0)}(W_{d_n}(\pi_{0,t_n}(y,\cdot),\pi_{0,t_n}(z,\cdot))   <\frac12\right]>\frac12.
$$
By Lemma \ref{l5.6i} and stationarity
$$
\begin{array}{l}
\ds \widetilde\P\left[ W_{d_n}(\mu_{0}^1, \mu_{0}^2) < 1  \right]\\
\\
=\widetilde\P\left[ W_{d_n}(\mu_{t_n}^1, \mu_{t_n}^2) < 1 \right]\\
\\
\ds\ge \widetilde\P\left[ W_{d_n}(\mu_{t_n}^1, \mu_{t_n}^2) < 1- \frac12 \mu_0^1(B(x_0,\gamma_0))
\wedge \mu_0^1(B(x_0,\gamma_0)) \right]\\
\\
 \ds>\frac12.
\end{array}
$$
Letting $n\to \infty$ yields
$$
\widetilde\P\left[\|\mu_{0}^1- \mu_{0}^2\|_{TV}<1 \right]>\frac12.
$$
Equivalently
$$
\lambda\{ h; \; |\mu_h^1- \mu_{h}^2\|_{TV}<1 \}>\frac12
$$
But this is impossible because $\mu^1_h$ and $\mu^2_h$ are almost surely singular. $\Box$

\section{Application to $2D$ Navier--Stokes equations}
 
We illustrate the above theory on the two-dimensional  Navier--Stokes equations on a bounded domain
$\mathcal O\subset \R^2$ with Dirichlet boundary conditions and a stationary forcing term.
The unknowns are the velocity $X(t,\xi)$ and the pressure
 $p(t,\xi)$ defined for $t>0$ and $\xi\in \overline{\mathcal O}$:
\begin{equation}
\label{e3.1} 
\left\{\begin{array}{ll}
 dX(t,\xi)=&[\Delta
X(t,\xi)-(X(t,\xi)\cdot\nabla)X(t,\xi)+g(X(t,\xi),\overline Y(t,\xi))]dt\\
\\
&-\nabla
p(t,\xi)dt+f(t,\xi)dt+ \sigma(X(t,\xi),\overline Y(t,\xi))dW,
\\
\mbox{\rm div}\;X(t,\xi)=&0,
\end{array}\right.
\end{equation}
with Dirichlet boundary conditions
$$
X(t,\xi)=0,\quad t>0,\;\xi\in \partial \mathcal O,
$$
and supplemented with the initial condition
$$
X(s,\xi)=x(\xi),\;\xi\in  \mathcal O.
$$

For simplicity, we assume that $f=0$. Following the usual notations we rewrite the equations as,
\begin{equation}
\label{e3.2}
\left\{\begin{array}{lll}
dX(t)&=&(AX(t)+b(X(t))+g(X(t),\overline Y(t)))dt\\
\\
&&+\sigma(X(t),\overline Y(t))dW(t),\quad s\le t,\\
\\
X(s)&=&x.
\end{array}\right.
\end{equation}
Here $A$  is the Stokes operator
 $$A=P\Delta,\quad D(A)=(H^2(\mathcal O))^2\cap (H^1_0(\mathcal O))^2\cap H,$$
where
$$
H=\{x\in (L^2(\mathcal O))^2:\;
\mbox{\rm div}\;x=0\;\mbox{\rm in}\;\mathcal O \},
$$
$P$ is the orthogonal projection of $(L^2(\mathcal O))^2$ on $H$ and     $b$  the operator
$$
(b(y),z)=b(y,y,z),\quad y,z\in  
V=\{y\in (H^1_0(D))^2 \cap H:\;\nabla\cdot y=0  \},
$$
where
$$
b(y,\theta, z)=-\sum_{i,j=1}^{2}\int_D y_i\;D_i\theta_j\;z_j\;d\xi,\quad y,\theta, z\in V.
$$
Moreover $W$ is a cylindrical Wiener process on a filtered probability space
$(\Omega,\mathcal F, (\mathcal F_t)_{t\ge 0},\P)$
in $H.$     
Finally,
$f\colon \R\to H$ is continuous and $2\pi$--periodic with respect to $t$.
We denote by $|\cdot|$ the norm in $H$,
by   $\|\cdot\|$ the norm in $V$ and by $(\cdot,\cdot)$ the scalar product in $H$.\bigskip

We assume that there exist $K_1>0$ and $h\;: \;\R\to\R$ such that
\begin{equation}
\label{e3.3}
\begin{array}{l}
\ds|g(x,y)|_H\le h(|y|_K),\quad\forall\; x\in H,\;y\in K,\\
\\
\ds|\sigma(x,y)|_{\mathcal L_2(H)}\le h(|y|_K),\quad\forall\; x\in H,\;y\in K,\\
\\
\ds \widetilde\E\left( h(|\overline Y(t)|_K)^2\right)=K_1<\infty.
\end{array} 
\end{equation}
\begin{Proposition}
\label{p3.1}
There exists an invariant measure for $Z$.
\end{Proposition}
{\bf Proof}.
By It\^o's formula we deduce that
\begin{equation}
\label{e3.4}\ds
\E|X(t)|^2+\E\int_0^t\|X(s)\|^2ds\ds
\le |x|^2+2\E\int_0^th(|\overline{Y}(s)|)^2ds
\end{equation}
and, by stationarity of $\overline Y$ and \eqref{e3.3},
\begin{equation}
\label{e3.5}
\E_1|X(t)|^2+\E_1\int_0^t\|X(s)\|^2ds
\le |x|^2+2tK_1.
\end{equation}
Therefore
\begin{equation}
\label{e3.6}
\frac1t\;\E_1\int_0^t\|X(s)\|^2ds\le |x|^2+2K_1.
\end{equation}
We deduce from \eqref{e3.6}, for any $M>0$,
\begin{equation}
\label{e3.7}
\frac1t\;\int_0^t\P_1(X(0,-s;0)\in B_{H^1}(0,M))ds\ge 1-\frac{2K_1}{M^2}.
\end{equation}
On the other hand, for any $\epsilon>0$ there exists a compact set $A_\epsilon\subset \mathcal K)$ such that
$$
\P(\overline{Y}_{(-\infty,0]}\in A_\epsilon)\ge 1- \epsilon
$$
and by the stationarity of $\overline{Y}$,
$$
\P(H(\theta,-t;\overline{Y}_{(-\infty,0]}\in A_\epsilon)\ge 1- \epsilon,\quad\forall\;  t\in\R.
$$
By the Krylov--Bogoliubov theorem there exists an invariant measure $\nu$ for $Z$. $\Box$

\bigskip

By Proposition \ref{p3.1} we deduce the existence of the family
$(\widetilde\mu_t)_{t\in\R,\widetilde{\omega}\in \widetilde{\Omega}}$.  Concerning uniqueness, we can use the criteria derived in sections \ref{s4} and \ref{s5}. For instance, assume that the noise 
is addiditive : $\sigma(x,y)=\sqrt C$ and non degenerate in the following sense
\begin{equation}
\label{e5.2g}
\mbox{\rm\;Tr}\; C<\infty,\quad C^{-1/2}(-A)^{-1/2}\in L(H).
\end{equation}
Then using the arguments in \cite{DPD-M2AN-control}, \cite{DPD-NS3D}, \cite{DO-JEE}, \cite{maslowski-flandoli} it is not difficult to 
prove that the strong Feller property holds and that the transition semigroup is almost surely irreducible. Then by Theorem \ref{t4.3g} and section \ref{s3}, there exists a unique 
evolutionary system of measure.

The nondegeneracy assumptions can be weakened using the ASF property. If we consider 
periodic boundary condition instead of Dirichlet boundary then the argument in \cite{HM} section 4.5 can be adapted to our setting and prove that the ASF property holds if 
$$
\E\left(\exp\left( th(\overline Y(t))^2\right)\right)\le c_1 \exp \left(c_2 t +1\right)
$$
for some $c_1,c_2\ge 0$ and if the noise acts on a sufficiently large number of modes.
Then, uniqueness of evolutionary system of measures holds if one can prove that there
exists $x_0$ such that \eqref{support} holds. 

It is probably also possible to extends the more difficult truly elliptic  case treated in  \cite{HM} section 4.6 but this requires much more work and is beyond the scope of this work.

\subsection{Uniqueness by coupling}

In the same spirit as in \cite{DD}, we show that coupling arguments extend to our situation. 
We do not consider the most general case which requires lengthy proofs. Instead, we consider a non degenerate noise so that the argument is not too long. The case of degenerate noise treated for instance in \cite{mattingly}, \cite{kuksinBook}, \cite{kuksin2}, \cite{mattingly2} could be treated by mixing the  arguments in these papers 
and the ideas below.

We are here concerned with the equation \eqref{e3.2} with a non degenerate additive noise
$\sigma(x,y)=\sqrt C$ satisfying \eqref{e5.2g}.

We need a further assumption on the process $\overline Y$. We assume that it posseses a Lyapunov 
structure. More precisely that there exists $\kappa_1,\; \kappa_2>0$ such that for all 
$t\ge s$:
\begin{equation}
\label{e7.9}
\widetilde\E\left( |\overline Y (t)|^2|{\mathcal F_s}\right)\le 
e^{-\kappa_1(t-s)}|\overline Y (s)|^2+\kappa_2.
\end{equation}
Also, for simplicity, we consider the case when \eqref{e3.3} holds with $h(r)=\kappa_3(1+x),\; x\ge 0$.
A different $h$ would require a different Lyapunov structure.

By Ito's formula and Poincar\'e inequality:
$$
\E(|X(t,s,x)|^2)\le e^{-\lambda_1 (t-s) } |x|^2 + \frac{2\kappa_3}{\lambda_1} + 2\kappa_3\int_s^t
e^{-\lambda_1 (t-\sigma) } |\overline Y(\sigma)|^2 d\sigma,
$$
where $\lambda_1$ is the first eigenvalue of $A$. Assumption \eqref{e7.9} implies:
$$
2\kappa_3\widetilde\E\left( \int_s^t e^{-\lambda_1(t-s)}|\overline Y (\sigma)|^2d\sigma |{\mathcal F_s}\right)\le 
\frac{2\kappa_3}{\kappa_1+\lambda_1} |\overline Y (s)|^2+\frac{2\kappa_2\kappa_3}{\lambda_1}.
$$
By the Markov property, we obtain 
for $t\ge r\ge s$
$$
\begin{array}{l}
\ds \E_1\left( |X(t,s,x)|^2+\delta  |\overline Y(t)|^2|\mathcal F^1_r\right)\\
\\
\ds \le e^{-\lambda_1(t-r)} |X(r,s,x)|^2 +\frac{2\kappa_3}{\lambda_1} + \frac{2\kappa_3}{\kappa_1+\lambda_1} |\overline Y (s)|^2+\frac{2\kappa_2\kappa_3}{\lambda_1}+\delta e^{-\kappa_1(t-r)}|\overline Y(r)|^2  +\delta\kappa_2.
\end{array}
$$ 
Recall that $\mathcal F^1_r=\widetilde{\mathcal F}_r\times \mathcal F_r$.
Let $T\ge 0$ and set 
$$
\kappa_5=\min\{\lambda_1,\frac12\kappa_1\},\;\alpha =e^{-\kappa_5 T}-e^{-\kappa_1 T},\;\delta = \frac{2\kappa_3}{\alpha(\kappa_1+\lambda_1)},\;
\kappa_4= \frac{2\kappa_3}{\lambda_1} +\frac{2\kappa_2\kappa_3}{\lambda_1}+\delta\kappa_2
$$
then, for $k\ge 0$, we obtain
$$
\begin{array}{l}
\ds\E_1\left( |X((k+1)T+s,s,x)|^2+\delta  |\overline Y((k+1)T+s)|^2|\mathcal F^1_{kT+s}\right)\\
\\
\ds\le e^{-\kappa_5 T}\left(|X(kT+s,s,x)|^2 + \delta  |\overline Y(kT+s)|^2\right) +\kappa_4.
\end{array}
$$
Thanks to this inequality, Lemma \ref{l5.1g}  below can now be proved as \cite[Lemma 6.1]{DD}.
\begin{Lemma}
\label{l5.1g}
Fix $T>0, M\in \N$ and set
$$
\tau=\inf\{kT+s:\;k\in \N, \;|X(kT+s,s,x)|^2+ \delta  |\overline Y(kT+s)|^2\le M\kappa_4\}.
$$
Then there exists $C(T)$, $M(T)\in \N$ such that if $M\ge M(T)$, we have
\begin{equation}
\label{e5.3g}
\P_1(\tau\ge kT+s)\le C(T) e^{-\frac12\;k\kappa_5 T}(1+|x|^2)
\end{equation}
and there exists $C(\alpha, T)$ such that for $\alpha<\frac12\;\kappa_5$
\begin{equation}
\label{e5.4g}
\E_1(e^{\alpha\tau})\le C(\alpha, T)e^{\alpha s}(1+|x|^2).
\end{equation}
 \end{Lemma}
The following Lemma is also proved as in \cite[Lemma 6.2]{DD}.
\begin{Lemma}
\label{l5.2g}
For any $x,y\in H,\;t\ge 0,\;s\in \R$, $K>0,\;g\in C^1_b(H)$ such that  $\|g\|_0\le 1$
there exist $c_1,c_2>0$ such that
\begin{equation}
\label{e5.5g}
\begin{array}{l}
|P_{s,t+s}g(x)-P_{s,t+s}g(y)|\le\ds
\frac1K\;e^{-\lambda_1t}(|x|^2+|y|^2)+\frac{\mbox{Tr } C}{K}\\
\\
\ds+\frac2{\lambda_1K}\;\int_s^{s+t}e^{-\lambda_1(t+s-\sigma)}|\overline Y(\sigma)|^2d\sigma +\frac{c_2}t\;e^{2c_1K} |x-y|.
\end{array}
\end{equation}
\end{Lemma}
\begin{Corollary}
\label{c5.3g}
For any $t>0$ there exists $\delta>0$ such that for any $s\in \R$,
\begin{equation}
\label{e5.6g}
\widetilde \E \left(\|P_{s,t+s}^*\delta_x-P_{s,t+s}^*\delta_y\|_{TV}   \right)\le \frac12,
\end{equation}
provided $|x|,|y|\le \delta.$
\end{Corollary}
{\bf Proof}. Take the supremum in $\|g\|_0\le 1$ and then the expectation $\widetilde \E$ in \eqref{e5.5g}. The result
follows taking first $K$ large and then choosing $\delta$ small.  $\Box$

\begin{Lemma}
\label{l5.4}
For any $T>0,\;\rho_1>0,\;\delta_1>0$ there exists
$K_0(\rho_1,\delta_1)\in \N$ and $\alpha(\rho_1,\delta_1)>0$
such that for any $|x|\le \rho_1$
\begin{equation}
\label{e5.7g}
\P_1(|X(K_0(\rho_1,\delta_1)T+s,s,x)|\le \delta_1)\ge \alpha(\rho_1,\delta_1).
\end{equation}
\end{Lemma}
{\bf Proof}. The proof is similar to that of \cite[Lemma 7.4]{DD}. The left hand side in \eqref{e5.7g} is independent on $s$. It suffices to consider 
$s=0$. We consider the deterministic problem
\begin{equation}
\label{e5.8g}
\left\{\begin{array}{l}
\ds\frac{d\widetilde{X}}{dt}=A\widetilde{X}+b(\widetilde{X}),\\
\\
\widetilde{X}(0)=x.
\end{array}\right.
\end{equation}
It is easy to see that for any $\rho_1>0, \delta_1>0$ there exists $K_0(\rho_1 , \delta_1)\in \N$
such that
\begin{equation}
\label{e5.9g}
|\widetilde{X}(t,0,x)|^2\le\frac14\;\delta_1^2\quad\mbox{\rm for}\;t\ge K_0(\rho_1 , \delta_1)T,
\end{equation}
provided $|x|\le \rho_1$.
Let
\begin{equation}
\label{e5.10g}
R(t)=\int_0^{t}e^{(t-\sigma)A}\overline Y(\sigma)d\sigma+\int_0^{t}e^{(t-\sigma)A}\sqrt{C}\;dW(\sigma).
\end{equation}
Since the noise is non degenerate we can prove that
$$
\P\left(\sup_{t\in K_0(\rho_1,K_1)T}|R(t)|\le \frac\eta2   \right)\ge \gamma_1(\eta,\widetilde\omega)>0,
\; \widetilde \P\mbox{-a.s}.
$$
Taking the expectation $\widetilde \E$ we find
$$
\P_1\left(\sup_{t\in K_0(\rho_1,K_1)T}|R(t)|\le \frac\eta2   \right)\ge \widetilde \E(\gamma_1(\eta,\widetilde\omega))>0.
$$
Now we can conclude the proof as in \cite{DD} with minor modifications. $\Box$

We are now ready to construct a coupling by proceeding as in the proof of \cite[Proposition 6.5]{DD}.
\begin{Proposition}
\label{p5.5g}
There exist $c>0$ and $\gamma>0$ such that for any $s\in \R,k\in \N$
and any $\varphi\in C_b(H)$
\begin{equation}
\label{e5.11g}
 \widetilde\E(|P_{s,t}\varphi(x)-P_{s,t}\varphi(y)|) \le c\|\varphi\|_0e^{-\widetilde\gamma (t-s)}(1+|x|^2+|y|^2).
\end{equation}
\end{Proposition}
{\bf Proof}. 
Fix $T>0$. We  take $\delta>0$ as in Corollary \ref{c5.3g}. For   $x,y\in B_\delta$ we fix $\widetilde\omega\in \widetilde\Omega$ and we choose a maximal coupling,  $(Z_1^s(x,y),Z_2^s(x,y))$. Then
$$
\P(Z_1^s(x,y)\neq Z_2^s(x,y))=\|P^*_{s,t}\delta_x-
P^*_{s,t}\delta_t\|_{TV}.
$$
By Corollary \ref{c5.3g} it follows that
$$
\P_1(Z_1^s(x,y)\neq Z_2^s(x,y))\le \frac12.
$$
We notice that $(Z_1^s(x,y),Z_2^s(x,y))$ is still a coupling of $(X(T+s,s,x),X(T+s,s,y))$ considered as random variables on  $\omega_1=(\omega,\widetilde\omega)$. Now we continue as in \cite[Proposition 6.5]{DD}, proving finally that
$$
\P_1(X_1^{s,h}\neq X_2^{s,h})\le ce^{-\widetilde\gamma k}(1+|x|^2+|y|^2),
$$
for some $\widetilde \gamma$. It follows:
$$
 \widetilde\E(|P_{s,kT+s}\varphi(x)-P_{s,kT+s}\varphi(y)|) \le c\|\varphi\|_0e^{-\widetilde\gamma k}(1+|x|^2+|y|^2).
$$
Writing
$$
P_{s,t}\varphi(x)-P_{s,t}\varphi(y)=P_{s,kT+s}P_{kT+s,t}\varphi(x)-P_{s,kT+s}P_{kT+s,t}\varphi(y)
$$
and $\|P_{kT+s,t}\varphi\|_0\le \|\varphi_0\|_0$, \eqref{e5.11g}   follows with a different constant $c$. $\Box$

By Borel-Cantelli Lemma, we deduce $\widetilde \P$ almost sure exponential convergence.
\begin{Corollary}
For $\widetilde \P$ almost  every $\widetilde \omega$, there exists $T_0(\widetilde\omega)$ such that
for any $t\ge T_0$
$$
|P_{s,t+s}\varphi(x)-P_{s,t+s}\varphi(y)| \le c\|\varphi\|_0e^{-\widetilde\gamma t/2}(1+|x|^2+|y|^2).
$$
\end{Corollary}

\begin{Theorem}
Let 
$$
P_{s,t}\varphi(x)=\E[\varphi(X(t,s,x))],\quad \varphi\in B_b(H),
$$
where $X(t,s,x)$ is the solution of the Navier-Stokes (NS).
Then, if \eqref{e7.9} holds,  we have
\begin{equation}
\label{e4.10ebis}
{\widetilde \E}\left|P_{s,t}\varphi(x)-\int_H\varphi(y)\mu_t(dy)\right|\le c(s)\|\varphi\|_0e^{-\widetilde\gamma (t-s)}(1+|x|^2)
\end{equation}
and
\begin{equation}
\label{e4.10e}
\left| P_{s,t}\varphi(x)-\int_H\varphi(x)\mu_{t}(dx)\right|\le K(\widetilde\omega,s,x)\|\varphi\|_0 e^{-\widetilde\gamma t/2}\quad\widetilde\P\mbox{\rm -a.s.}.
\end{equation}
\end{Theorem}
{\bf Proof}.  Taking into account \eqref{e3.1f} we have
$$
P_{s,t}\varphi(x)-\int_H\varphi(y)\mu_t(dy)=
\int_H(P_{s,t}\varphi(x)-P_{s,t}\varphi(y))\mu_s(dy).
$$
Therefore by \eqref{e5.11g} it follows that
$$
{\widetilde \E}\left|P_{s,t}\varphi(x)-\int_H\varphi(y)\mu_t(dy)\right|\le c\|\varphi\|_0e^{-\widetilde\gamma (t-s)}
\int_H(1+|x|^2+|y|^2)\mu_s(dy)
$$
so that \eqref{e4.10ebis} follows.  Again, \eqref{e4.10e} is obtained thanks to Borel Cantelli Lemma. $\Box$

Identity \eqref{e4.10e} can be interpreted by saying that as $t\to+\infty$, the observable $ P^{\widetilde\omega}_{s,t}\varphi$ approaches exponentially fast a random limit curve
$$
  t\to  \int_H\varphi(x)\mu_t(dx), 
$$
which forgets about $s$.

\end {document}